\documentclass{article}

\usepackage{epsfig}
\usepackage{graphicx}
\makeatletter
\@ifundefined{showcaptionsetup}{}{%
 \PassOptionsToPackage{caption=false}{subfig}}
\usepackage{subfig}
\makeatother

\usepackage{amssymb}
\usepackage{amsthm}
\usepackage{amsmath}
\usepackage{hyperref}
\newtheorem{thm}{Theorem}

\title{Localized Harmonic Characteristic Basis Functions for Multiscale Finite Element Methods}

\author{Leonardo A. Poveda\thanks{Instituto de Matem\'atica e Estat\'istica, Universidade de S\~ao Paulo, Brazil, email: \href{mailto:lpovedac@ime.usp.br}{\texttt{lpovedac@ime.usp.br}} } 
\and Sebastian Huepo\thanks{Departamento de Matem\'aticas, Universidad Nacional de Colombia, Bogot\'a D.C., Colombia,  email: \href{mailto:shuepobe@unal.edu.co}{\texttt{shuepobe@unal.edu.co}} }
\and Juan Galvis\thanks{Departamento de Matem\'aticas, Universidad Nacional de Colombia, Bogot\'a D.C., Colombia, email: \href{mailto:jcgalvisa@unal.edu.co}{\texttt{jcgalvisa@unal.edu.co}} }
\and Victor M. Calo \thanks{Center for Numerical Porous Media, Applied Mathematics.  Computational Science and Earth Sciences \& Engineering, King Abdullah University of Science and Technology. Thuwal 23955-6900 Kingdom of Saudi Arabia, \href{mailto:victor.calo@kaust.edu.sa}{\texttt{victor.calo@kaust.edu.sa}} }
}

\begin{document}

\maketitle

\tableofcontents

\section*{Abstract}

We solve elliptic systems of equations posed on highly heterogeneous materials. Examples of this class of problems are  composite structures and geological processes.  
We focus on a  model problem which is a second-order elliptic equation with discontinuous coefficients. These coefficients represent the conductivity of a composite material. 
We assume a background with low  conductivity that contains inclusions with 
different thermal properties.
Under this scenario we design a multiscale finite element method to efficiently approximate solutions. 
The method is based on an asymptotic expansions of the solution in terms of the ratio between the conductivities. 
The resulting method constructs (locally)  finite element basis functions (one for each inclusion). These bases that generate the multiscale finite element space 
where the approximation of the solution is computed. Numerical experiments show the good performance of the
proposed methodology. \\

{\bf keywords:} Elliptic equation, asymptotic expansions, high-contrast coefficients, multiscale finite element method, harmonic characteristic function.

\section{Introduction}
\label{Sec1:Introduction}

Many physical and engineering applications naturally require multiscale solutions.  This is specially true for problems related to metamaterials, composite materials, and porous media flows; see \cite{chen2013multiscale,li2011finite,ozgun2013software,epov2015effective,zhou2012interaction}. The mathematical and numerical analyses for these problems are  challenging since they are governed by elliptic equations  with high-contrast coefficients 
(\cite{hou1997multiscale,chen2003mixed,galvis2010domain,calo2014asymptotic}). 
For instance, in the modeling of composite materials, their conducting or elastic properties are modeled by discontinuous coefficients. The value of the coefficient can vary several orders of magnitude across discontinuities. Problems with these jumps are referred to as high-contrast. Similarly, the coefficient is denoted as a high-contrast coefficient. See for instance 
 \cite{calo2014asymptotic,galvis2014spec,efendiev2013generalized, efendiev2012coarse}.

We seek to understand how the high-contrast variations in the material properties affect the structure of the solution. In terms of the model, these variations appear in the coefficients of the differential equations. We expand our previous work \cite{calo2014asymptotic}, where we construct an asymptotic expansion to represent solutions. The asymptotic expansion is obtained in terms of  the high contrast in the material properties.
In this paper, we use this asymptotic expansion design numerical solutions for high-contrast problems. The asymptotic expansion helps us derive elegant numerical strategies and to understand the local behavior of the solution. In addition, the asymptotic expansion can be used to study functionals of solutions and describe their behavior with respect to the contrast or other important parameters.
The asymptotic expansion  in  \cite{calo2014asymptotic} uses 
globally supported harmonic extensions of subdomains indicator functions referred to as \emph{harmonic characteristic functions}.
 We modify the construction presented in \cite{calo2014asymptotic} to approximate the harmonic characteristic function in the local neighborhood of
each inclusion. Thus, in order to make practical use of the expansion we avoid computing each characteristic function for the whole domain. This modification renders the method computationally tractable while the reduction in accuracy is not significant. For the case of dense distributions of inclusions, we observe numerically that the optimal size for the support of the basis functions is of the order of the representative distance between inclusions.
We perform numerical tests that show the good performance of the proposed  Multiscale Finite Element Method.
We use a finite element method (FEM) and assume that there is a fine-mesh that completely resolves 
the geometrical configuration of the inclusions in the domain. That is, the fine-scale finite element formulation fully captures the solution behavior. To compute the linear system solution at this fine resolution is not practical and therefore, a multiscale finite element strategy is needed in order to compute a coarse-scale representation that captures relevant information of the targeted fine-scale solution. The coarse dimension in our simulations corresponds to the total number of inclusions.
Nevertheless, a  coarser scale may be needed for some applications. In this case it is possible to use the framework of the generalized finite element method to design and analyze a coarser scale for computations.
For a detailed discussion, see \cite{efendiev2013generalized} and references therein.
In some more demanding applications an efficient iterative 
domain decomposition method could also be designed and analyzed for these problems.  This is under investigation 
and will be presented elsewhere.

The problem of computing solutions of elliptic problems related to modern artificial materials such as dispersed and/or densely packed composite materials
has been considered by some researchers recently. For instance in  \cite{wen2004effective} the computation of effective properties of dispersed composite materials is considered. They use a classical  multiscale finite element method. We recall that the application of the classical finite element method may lead to the precense of \emph{resonance} errors 
due to the chosen local boundary conditions. See  \cite{MR2830328}.
In \cite{p2013, peter1} the authors develop a finite element method based on a network approximation of the conductivity for particle composites.

The rest of the paper is organized as follows. In Section 2, we setup the problem. 
In Section 3, we summarize  the asymptotic expansion procedure described in \cite{calo2014asymptotic}. We also introduce the definition of
 harmonic characteristic functions, which help us determine the individual  terms of the expansion. 
 In Section 4, we illustrate some aspects of the asymptotic expansion using some finite element computations.
 Section 5 constructs multiscale finite elements using the asymptotic expansion described in the previous sections. In particular, we approximate the leading term of the expansion with localized harmonic characteristic functions. 
We then apply this approximation to the case of dense high-contrast inclusions. In Section 6, we present some 
numerical experiments using the methods proposed. Finally, in Section 7 we draw  some conclusions.

\section{Problem Setup}\label{sec2:problem setting}

We use the notation introduced in \cite{calo2014asymptotic}. We consider a second order  elliptic problems of the form,
\begin{equation} \label{Strongform}
-\mbox{div }(\kappa \nabla u)=f,\ \mbox{  in }  D,  
\end{equation}
with Dirichlet data defined by $u=g$ on $\partial D$. Here $D$ is the disjoint union of a background domain  $D_0$ and subdomains 
that represent the inclusions, i.e., $D=D_0\cup (\bigcup_{m=1}^M \overline{D}_m)$.
 We assume that $D_1,\dots,D_M$ are connected polygonal domains (or domains with smooth boundaries). Additionally, we  require  that each inclusion $D_m$, 
$m=1,\dots,M$ is compactly included in the open set $D\setminus 
\bigcup_{\ell=1, \ell \not=m}^M\overline{D}_\ell$, i.e., 
$\overline{D}_m\subset D\setminus \bigcup_{\ell=1, \ell\not=m}^M \overline{D}_\ell$. Let $\kappa$ be defined by 
\begin{equation}\label{eq:coeff1inc-multiple}
\kappa(x)=\left\{\begin{array}{cc} 
\eta,& x\in   D_m, ~~m=1,\dots,M, \\
1,& x\in  D_0=D\setminus \bigcup_{m=1}^M\overline{D}_m.
\end{array}\right.
\end{equation}
Following  \cite{calo2014asymptotic}, we represent the solution by an asymptotic 
expansion in terms of the contrast $\eta$. The expansion reads, 
\begin{equation}\label{expansionu_eta}
u_\eta=u_0+\frac{1}{\eta}u_1+\frac{1}{\eta^2}u_2+\dots=
\sum_{j=0}^\infty \eta^{-j} u_{j},
\end{equation}
 with coefficients $\{u_j\}_{j=0}^\infty\subset 
H^1(D)$  and such that they satisfy the following Dirichlet boundary 
conditions, 
\begin{equation}\label{Boundarycondi}
u_0=g \mbox{ on } \partial D \quad \mbox{ and }
\quad u_{j}=0 \mbox{ on } \partial D \mbox{ for } j\geq 1.
\end{equation}

\section{Asymptotic Expansion}\label{sec:expansion}

Now we compute the terms in the  asymptotic expansion (\ref{expansionu_eta}). 
For more details on the construction and related expansions we refer to \cite{calo2014asymptotic, poveda}.

First, we introduce the \emph{harmonic characteristic functions}.
Let $\delta_{m\ell}$ represent the Kronecker delta, which is equal to $1$ when $m=\ell$ and $0$ otherwise.
For each $m=1,\dots, M$ we introduce the \emph{harmonic characteristic function} of $D_m$, $\chi_{D_{m}}\in H_{0}^{1}(D)$, with the condition
\begin{equation}\label{chimulti}
\chi_{D_{m}}\equiv\delta_{m\ell} \mbox{ in }D_{\ell},\mbox{ for }\ell=1,\dots,M,
\end{equation}
and which is equal to the harmonic extension of its boundary data  to the interior of $D_{0}$. We then have,

\begin{align}
\Delta \chi_m&=0,  \mbox{ in  } D_0,\\ \nonumber
\chi_m&=0,  \mbox{ on }\partial D \mbox{ and } \partial D_\ell, \, m\neq \ell, \, \ell=1,2,\dots,M\\ \nonumber
\chi_m&=1,  \mbox{ on }\partial D_m. \label{Chimultiinclusions}
\end{align}

The function $u_0$ in (\ref{expansionu_eta}) can be explicitly written in term of 
the harmonic characteristic functions and a boundary corrector; see \cite{calo2014asymptotic}. 
In fact, we can write,
\begin{equation}\label{formulau0}
u_0=u_{0,0}+\sum_{m=1}^{M}c_m(u_0)\chi_{D_m},
\end{equation}
where $u_{0,0}\in H^1(D)$,  for $m=1,\dots,M$, and $u_{0,0}$ solves the following problem posed
in the background $D_0$,
\begin{align}
-\Delta u_{0,0}&=f, \mbox{ in } D_0,\\\nonumber
u_{0,0}&=g,  \mbox{ on } \partial D\\  
u_{0,0}&=0, \mbox{ on  }  \partial D_m, m=1,\dots,M.\nonumber
\end{align}\label{Def:u00}
$u_{0,0}$ is globally supported but it is forced to vanish in all interior inclusions boundaries. 
The constants  in (\ref{formulau0}) solve an $M$ dimensional linear system. Let 
 $\mathbf{c}=(c_1(u_0),\dots,c_M(u_0))\in \mathbb{R}^M$, then we have that  
\begin{equation}\label{eq:matrixformu0}
\mathbf{A}_{geom}\mathbf{c}=\mathbf{b},
\end{equation}
where $\mathbf{A}_{geom}=[a_{m\ell}]$ and $\mathbf{b}=(b_1,\dots,b_M)\in \mathbb{R}^M$ are defined by
\begin{eqnarray}\label{Def:aml}
a_{m\ell}&=&\int_{D}\nabla \chi_{D_m}\cdot \nabla \chi_{D_{\ell}},
\end{eqnarray}
and 
\begin{eqnarray}
b_{\ell}&=&\int_{D}f\chi_{D_{\ell}}-\int_{D_0}\nabla u_{0,0}\cdot \nabla\chi_{D_{\ell}},\label{Def:bl}
\end{eqnarray}
respectively. 
$\sum_{m=1}^{M}c_m(u_0)\chi_{D_m}$ is the Galerkin projection 
of $u_0-u_{0,0}$ into the space $\mbox{Span}\left\{\chi_{D_m}\right\}_{m=1}^{M}$.

Now for the sake of completeness, we briefly describe the next individual terms of the asymptotic expansion. 
We have for $j=1,2,\dots,$
\begin{equation}
u_{j}=\widetilde{u}_{j}+\sum_{m=1}^{M}c_{j,m}\chi_{D_{m}},\label{eq:formulauj}
\end{equation}
where  the function $\widetilde{u}_{j}$ is defined in three steps
\begin{enumerate}
\item  Solve a Neumann problem in each inclusion with data from $u_{j-1}$ with $j=1,2,\dots$. We describe the next terms in the asymptotic expansion. We have the restriction of $u_j$ to the subdomain $D_m$ with $m=1,\dots,M$, that is
\[
u_j=\widetilde{u}_j+c_{j,m},\quad\mbox{with }\int_{D_m}\widetilde{u}_j=0,
\]
and $\widetilde{u}_j$ satisfies the Neumann problem
\[
\int_{D_m}\nabla \widetilde{u}_j\cdot \nabla z=\int_{D_m}fz-\int_{\partial D_m}\nabla u_{j-1}^{(0)}\cdot n_{m}z,\mbox{ for all }z\in H^1(D_m),
\]
for $m=1,\dots,M$. From now on, we use the notation $w^{(0)}$, which means that the function $w$ is restricted to the domain $D_{0}$, that is $w^{(0)}=w|_{D_0}$. The constants $c_{j,m}$ will be chosen suitably.
\item  Solve a Dirichlet problem in the background  $D_0$ with data  $\widetilde{u}_{j}$  in each inclusion.
For $j=1,2,\dots$, we have that $u_j$ in $D_m$, $m=1,\dots,M$, then we find $u_j^{(0)}$ in $D_0$ by solving the Dirichlet problem
\begin{align}\label{Def:ujDiric}
\int_{D_0}\nabla u_{j}^{(0)}\cdot \nabla z&=0, &&\mbox{ for all }z\in H_0^1(D_0)\\\nonumber
u_j^{(0)}&=u_j\, (=\widetilde{u}_j+c_{j,m}), &&\mbox{on }\partial D_m,\, m=1,\dots,M,\\\nonumber
u_j^{(0)}&=0, &&\mbox{on }\partial D.\nonumber
\end{align}
Since $c_{j,m}$ are constants, we define their corresponding harmonic extension by $\sum_{m=1}^{M}c_{j,m}\chi_{D_m}$. So we rewrite
\begin{equation}\label{Def:uj}
u_j=\widetilde{u}_j+\sum_{m=1}^{M}c_{j,m}\chi_{D_m}.
\end{equation}
\item The $u_{j+1}$ in $D_m$ satisfy the following Neumann problem
\[
\int_{D_m}\nabla u_{j+1}\cdot \nabla z=-\int_{\partial D_m}\nabla u_j^{(0)}\cdot n_0z,\quad\mbox{for all }z\in H^1(D).
\]
The compatibility condition is satisfied for $\ell=1,\dots, M$, then
\begin{eqnarray*}
0  = \int_{\partial D_{\ell}}\nabla u_{j+1}\cdot n_{\ell} & = & -\int_{\partial D_{\ell}}\nabla u_j^{(0)}\cdot n_0\\
& = & -\int_{D_{\ell}}\nabla \left(\widetilde{u}_j^{(0)}+\sum_{m=1}^{M}c_{j,m}\chi_{D_m}^{(0)}\right)\cdot n_0\\
& = & -\int_{\partial D_{\ell}}\nabla \widetilde{u}_j^{(0)}\cdot n_0-\sum_{m=1}^{M}c_{j,m}\int_{\partial D_m}\nabla \chi_{D_m}^{(0)}\cdot n_0.
\end{eqnarray*}
\end{enumerate}
  A detailed description of the differential problems can be found in \cite{calo2014asymptotic}. The constants $\{c_{j,m}\}$ in (\ref{eq:formulauj}) are computed solving a linear system 
similar to the one defined above in (\ref{eq:matrixformu0}).  We have that $\mathbf{c}_j=(c_{j,1},\dots,c_{j,M})$ is the solution of the system
\[
\mathbf{A}_{geom}\mathbf{c}_j=\mathbf{y}_j,
\]
where
\[
\mathbf{y}_j=\left(-\int_{D_0}\nabla \widetilde{u}_{j}^{(0)}\cdot \nabla \chi_{D_1},\dots ,-\int_{D_0}\nabla \widetilde{u}_j^{(0)}\cdot \nabla \chi_{D_M}\right).
\]
In \cite{calo2014asymptotic, poveda} the authors prove that  the expansion (\ref{expansionu_eta}) converges absolutely in $H^1(D)$ for $\eta$ sufficiently large.

\begin{thm}\label{theorem}
Consider the problem (\ref{Strongform}) with coefficient (\ref{eq:coeff1inc-multiple}). The corresponding expansion (\ref{expansionu_eta}) with boundary condition  (\ref{Boundarycondi}) converges absolutely in $H^1(D)$ for $\eta$ sufficiently large. Moreover, there exist positive constants $C$ and $C_1$ such that for every $\eta>C$, we have
\[
\left\Vert u-\sum_{j=0}^{J}\eta^{-j}u_j\right\Vert_{H^1(D)}\leq C_1\left(\|f\|_{H^{-1}(D)}+\|g\|_{H^{1/2}(\partial D)}\right)\sum_{j=J+1}^{\infty}\left(\frac{C}{\eta}\right)^{j},
\]
for $J\geq 0$.
\end{thm}

\section{An Expansion}\label{sec:illustration}

In this section we illustrate  the expansion in two dimensions.
We use \textsc{MatLab} for the computations, see \cite{guide1998mathworks}. In particular, a few terms are computed numerically using the finite element 
method, see for instance \cite{johnson2012numerical,hughes2012finite}. In particular, we solve the sequence of problems 
posed in 
the background subdomain and in the inclusions. Our main goal is to device efficient numerical approximations for $u_0$ (and then for $u_\eta$ by using 
Theorem \ref{theorem}). 

\begin{figure}[!hbt]
%a
\begin{centering}
\subfloat[\label{a1}Geometry and mesh.]{
\begin{centering}
\includegraphics[width=0.38\textwidth]{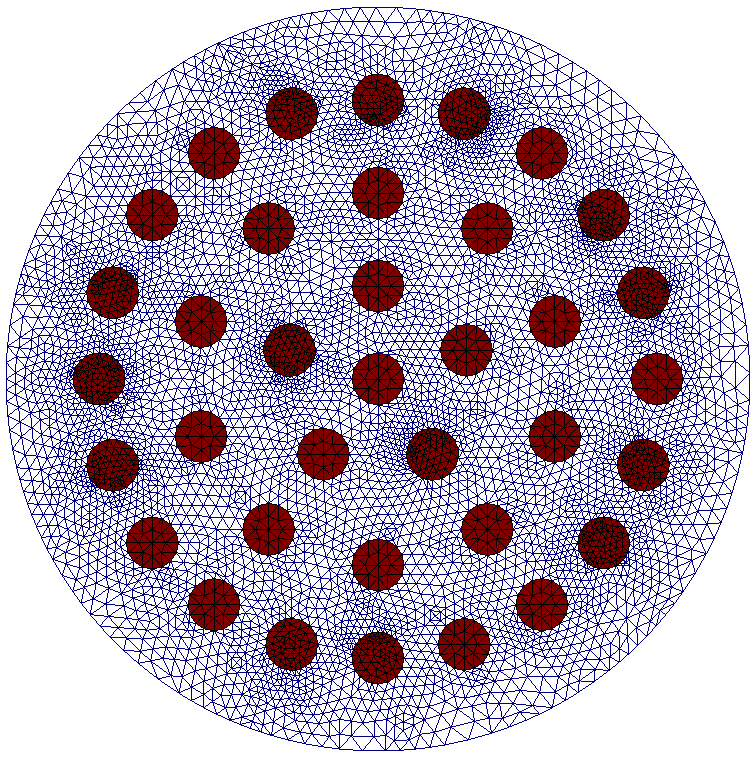}
\par\end{centering}
}
\par\end{centering}
%b
\begin{centering}
\subfloat[\label{b1}Finite element solution of problem \eqref{Example1Multiscale}  with  $\eta=6$.]{\begin{centering}
\includegraphics[width=0.50\textwidth]{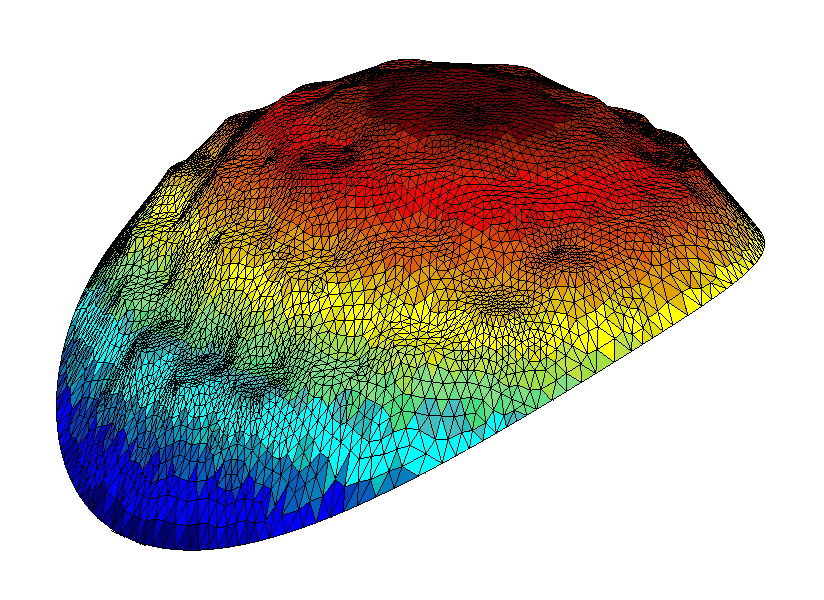}
\par\end{centering}
}
\par\end{centering}
%c
\begin{centering}
\subfloat[\label{c1} Asymptotic solution $u_0$  in \eqref{formulau0}]{\begin{centering}
\includegraphics[width=0.50\textwidth]{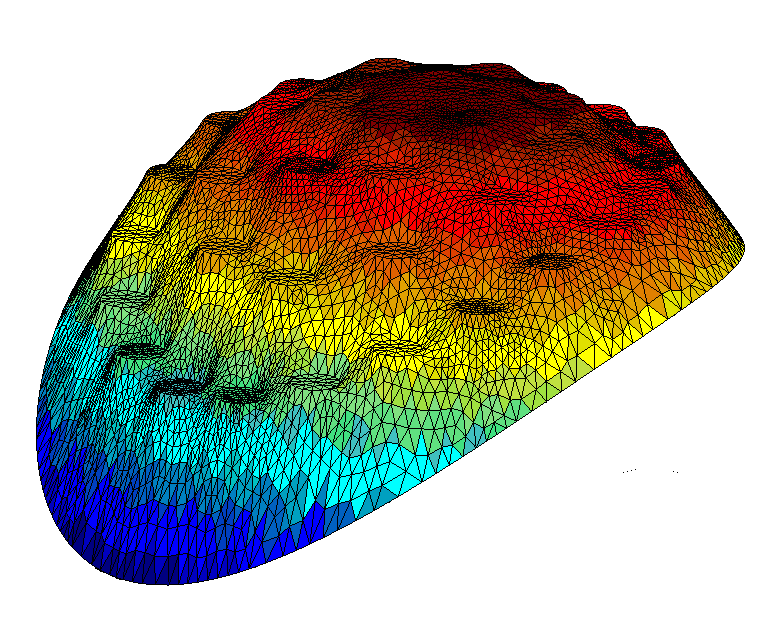}
\par\end{centering}
}
\par\end{centering}

\caption{Circular domain with $36$ identical inclusions. Geometry, direct numerical simulation with finite elements, and asymptotic expansion.}
\label{fig:figure1} 
\end{figure}

We consider $D=B(0,1)$ the circle with center $(0,0)$ and radius $1$.  We add $36$ (identical) circular inclusions of radius $0.07$. This is illustrated in the Figure \ref{a1}. Then, we numerically solve the problem
\begin{equation}\label{Example1Multiscale}
\left\{\begin{array}{ll}
-\mbox{div}(\kappa(x)\nabla u_\eta(x))=1, & \mbox{ in }D\\
\hspace{0.1in}u(x)=x_1+x_2^2, & \mbox{ on }\partial D.
\end{array}\right.
\end{equation}

Figure \ref{fig:figure1} shows that solution for $\eta=6$ against the computed $u_0$.

\begin{figure}[!htb]
\begin{centering}
\subfloat[\label{graphics-a}$\sum_{m=1}^{M}c_m(u_0)\chi_{D_m}$]{\begin{centering}
\includegraphics[scale=0.25]{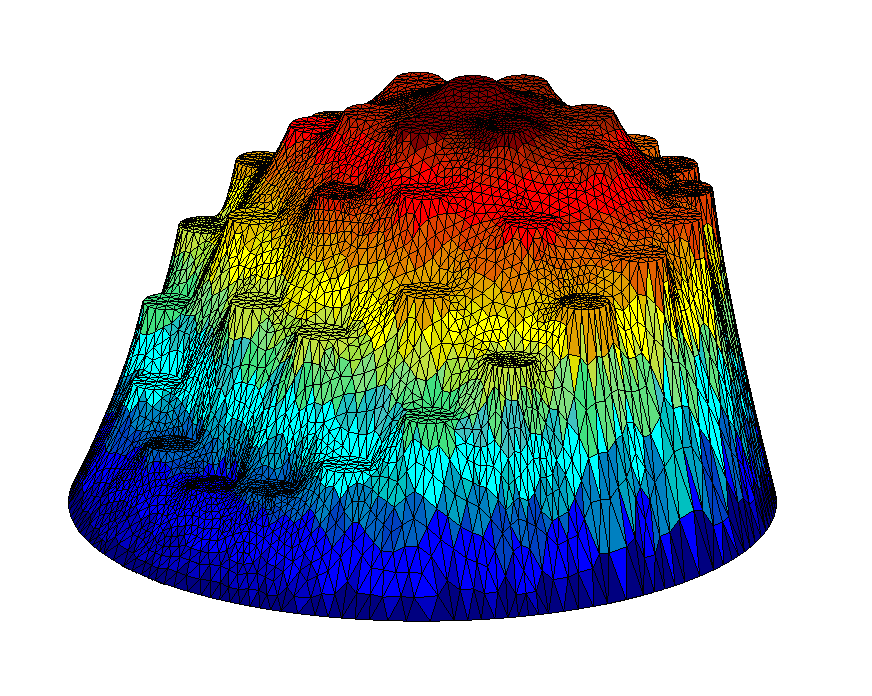}
\par\end{centering}
}
\par\end{centering}

\begin{centering}
\subfloat[\label{graphics-b}$u_{0,0}$]{\begin{centering}
\includegraphics[scale=0.25]{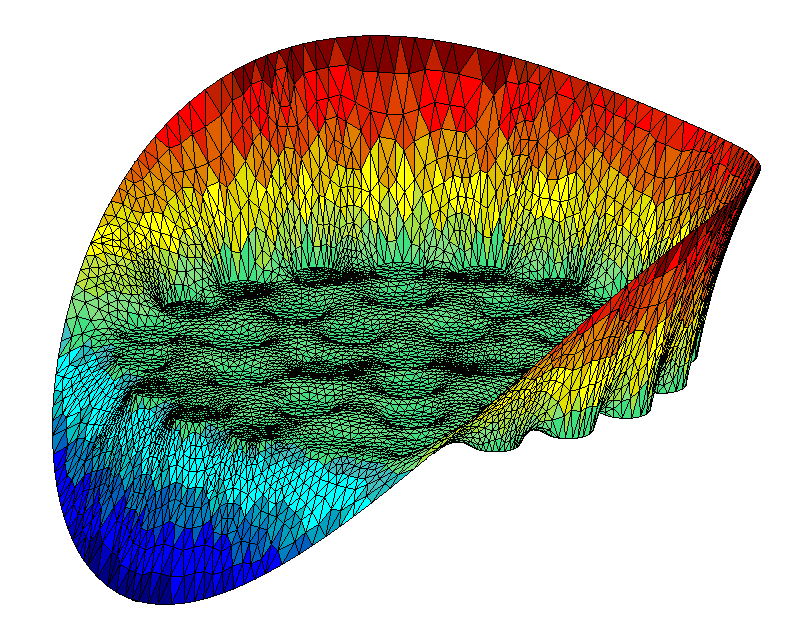}
\par\end{centering}
}
\par\end{centering}
\caption{Function $\sum_{m=1}^{M}c_m(u_0)\chi_{D_m}$ and $u_{0,0}$ for  problem (\ref{Example1Multiscale}) with $\eta=6$. See (\ref{formulau0}).}
\label{fig:figure1-1}
\end{figure}

In Figure \ref{fig:figure1-1} we show the two parts of $u_0$ in (\ref{formulau0}), the 
combination of the harmonic characteristic functions $\sum_{m=1}^{M}c_m(u_0)\chi_{D_m}$ and the boundary 
corrector $u_{0,0}$.  The results suggest that the boundary corrector $u_0$ decays fast away form the boundary 
$\partial D$.

\begin{figure}[!hbt]
\subfloat[$u_1$]{\includegraphics[width=0.42\textwidth]{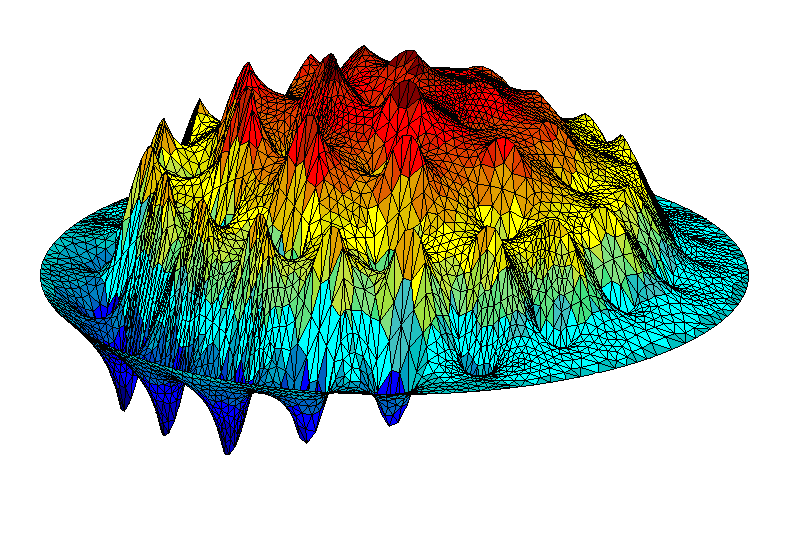}
}\subfloat[$u_2$]{\includegraphics[width=0.42\textwidth]{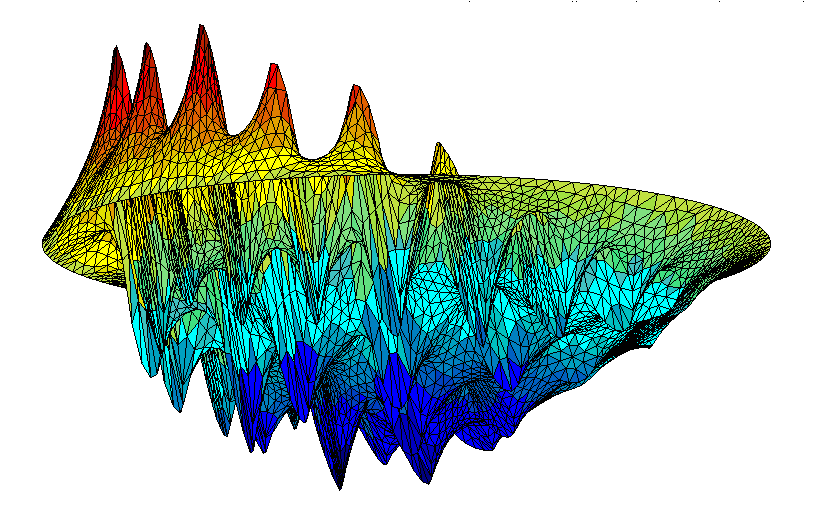}
}\\
\subfloat[$u_1$ restricted  to the inclusions]{\includegraphics[width=0.42\textwidth]{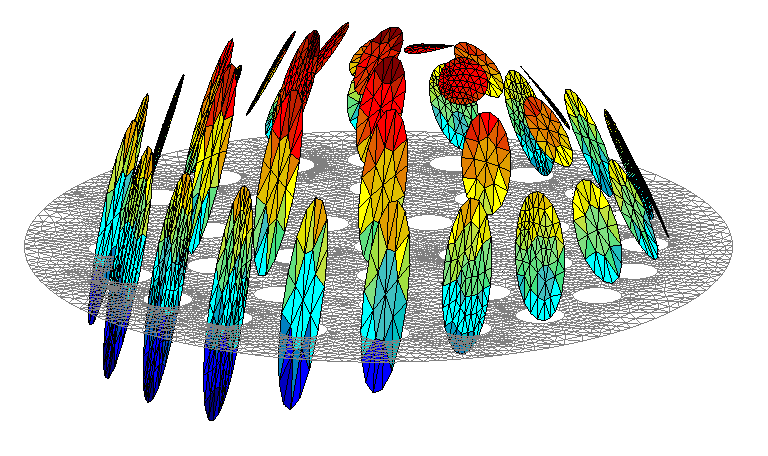}
}\subfloat[$u_2$ restricted to the inclusions]{\includegraphics[width=0.42\textwidth]{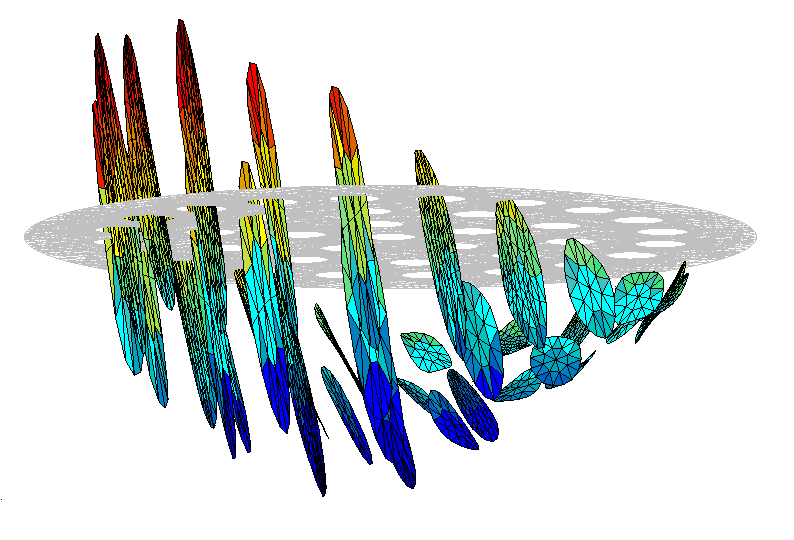}
}

\caption{Top: Functions $u_1$ and $u_2$ for the problem (\ref{Example1Multiscale}) with $\eta=6$. Bottom: 
Functions $u_1$ and $u_2$ restricted to the inclusions for the problem (\ref{Example1Multiscale}) with $\eta=6$.}
\label{fig:figure2}
\end{figure}

In Figure \ref{fig:figure2} we show the second and third term of the expansions,  $u_1$ and $u_2$. We also show the influence of $\eta$ on the convergence of the series  in Table \ref{table:3.1}. As predicted by Theorem \ref{theorem}, as $\eta$ grows convergence of the series expansion is accelerated. For example, problem (\ref{Example1Multiscale}) with $\eta=10$ requires eight terms in the series to achieve a relative error  of $10^{-8}$  while for $\eta=10^4$ requires only two terms.

\begin{table}[!h]\
\begin{center}\footnotesize
\begin{tabular}{|c|c|c|c|c|c|c|c|c|c|c|} \hline\hline
$\eta$ & 3 &6&10&$10^2$&$10^3$&$10^4$& $10^5$&$10^6$&$10^7$&$10^8$\\ \hline
\# & 25&    11&    8&     4&  
   3&     2&     2&     2&     1&     1 \\
\hline\hline
\end{tabular}
\end{center}
\caption{Number of terms needed to obtain a relative error 
of $10^{-8}$ for a given value of $\eta$.}\label{table:3.1}
\end{table}

\section{Approximation  with Localized Harmonic Characteristic Functions}\label{sec:approximation}

\begin{figure}[!htb]
\centering
\includegraphics[width=0.50\textwidth]{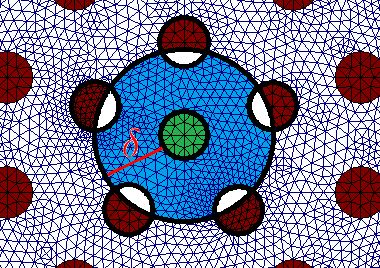}
\caption{Illustration of $\delta$-neighborhood of an inclusion. 
The selected inclusion is green.
The $\delta$-neighborhood of this inclusion is given in blue color, while preserving the truncated harmonic characteristics function.
We highlight with white color the other inclusions that are within the 
$\delta$-neighborhood of the selected inclusion.}
\label{fig:deltaneigh00} 
\end{figure}

In this section, we present a computable method based on the asymptotic expansion described in section \ref{sec:expansion} and the insights on the structure of the solution described in section \ref{sec:illustration}. We assume that $D$ is the union of a background and multiple inclusions that are homogeneously distributed. We only approximate of the leading term $u_0$. The
 remaining terms can be 
approximated similarly. We describe $u_0$ with localized harmonic characteristic functions. The computation of the
 harmonic characteristic functions is computationally expensive since these are fully global functions. That is, we approximate harmonic 
characteristics functions by solving a local problem (instead of a whole background problem). 
For instance, we pose a problem in a small neighborhood of each
inclusion.  The domain where the harmonic characteristic functions are computed is sketched in Figure 
\ref{fig:deltaneigh00}. The domain marked in Figure \ref{fig:deltaneigh00} corresponds to the adopted neighborhood
 of the inclusion painted with 
 green color. In this case, the approximated (or truncated) harmonic characteristic function is set to be
 zero on the boundary of the neighborhood of the selected inclusion.

The harmonic characteristic functions solve a background problem, which is global.  To reduce the cost of computing the characteristic harmonic functions, we solve problem (\ref{Example1Multiscale}), but restrict it to a neighborhood of the corresponding inclusion. The exact characteristic functions are defined by 
(\ref{chimulti}). We define the neighborhood of the inclusion $D_m$ by 
\[
D_{m,\delta}=\overline{D}_{m}\cup \left\{x\in D_0:d(x,D_m)<\delta\right\},
\]
and approximate the characteristic function for the $\delta$-neighborhood solving
\begin{align}
\Delta \chi_m^{\delta}&=0, \mbox{ in } D_{m,\delta},\nonumber \\
\chi^{\delta}_m&=0, \mbox{ on }\partial D_{m,\delta} \mbox{ and } \partial D_{\ell} \cap  D_{m,\delta} 
\mbox{ for } \ell \not = m, \nonumber \\ 
\chi^{\delta}_m&=1,  \mbox{ on }\partial D_m. \label{eq:deltachi}
\end{align}
The exact expression for $u_0$ is given by
\begin{equation}\label{u0multiscale}
u_0=u_{0,0}+\sum_{m=1}^{M}c_m\chi_{D_m}= 
u_{0,0}+u_c,
\end{equation}
where we have introduced $u_c=\sum_{m=1}^{M}c_m(u_0)\chi_{D_m}$. The matrix problem for 
$u_c$ with globally supported basis 
was given in (\ref{eq:matrixformu0}).
We now define $u_0^{\delta}$, the mulsticale approximation of $u_0$,  using a similar expression which is given by 
\begin{equation}\label{u0trnmultiscale}
u_0^{\delta}=u_{0,0}^{\delta}+
\sum_{m=1}^{M}c_m^{\delta}\chi_{D_m}^{\delta}=
u_{0,0}^{\delta}+u_c^{\delta},
\end{equation}
where $u_c^{\delta}=\sum_{m=1}^{M}c_m^{\delta}\chi_{D_m}^{\delta}$ with each
 $c_m^{\delta}$ computed similarly to $c_{m}$ using an alternative matrix problem with basis $\chi_{D_m}^{\delta}$ instead of $\chi_{D_m}$. This system is given by
\[
\mathbf{A}^{\delta}\mathbf{c}^{\delta}=\mathbf{b}^{\delta}.
\]
with $\mathbf{A}^{\delta}=[a_{m\ell}^{\delta}]$, with $a_{m\ell}^{\delta}=\int_{D_m}\nabla \chi^{\delta}_m
 \cdot \nabla \chi^{\delta}_{\ell}$, $\mathbf{c}=[c^{\delta}_0(u_0),\dots,c^{\delta}_M(u_0)]$ and
 $\mathbf{b}^{\delta}=[b_{\ell}^{\delta}]=\int_Df\chi_{D_{\ell}}^{\delta}$.

We also introduced the $u_{0,0}^{\delta}$, that is an approximation to the boundary corrector $u_{0,0}$, 
that solves the problem
\begin{align*}
-\Delta u_{0,0}^{\delta}& = f, \mbox{ in } D_0^{\delta},\\
u_{0,0}^{\delta}& =g,\mbox{ on }\partial D\\ 
u_{0,0}^{\delta}&=0,\mbox{ on }\partial D_0^{\delta}\cap D,
\end{align*}
where $D_0^{\delta}$ is the subdomain  within a distance   $\delta$  from the boundary $\partial D$, 
\[
D_0^{\delta}=\left\{x\in D:d(x,\partial D)<\delta\right\}.
\]

\section{Numerical Experiments}

\begin{figure}[!htb]
\centering
\subfloat[\label{a2} Illustration of the Global characteristic function for the problem \eqref{Example1Multiscale}.]{\includegraphics[width=0.4\textwidth]{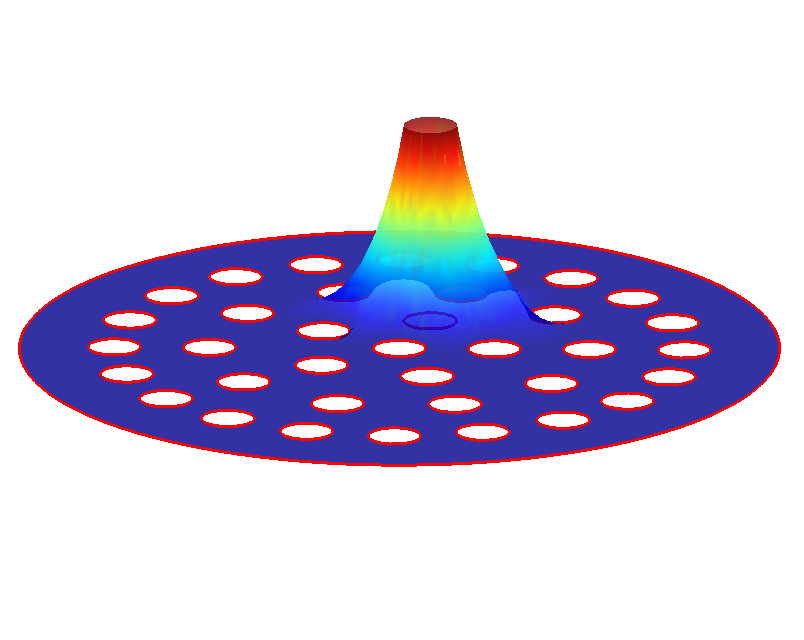}
}\hspace{1cm}\subfloat[\label{b2} Difference between the global characteristic function \eqref{chimulti} and the localized characteristic function \eqref{eq:deltachi}.]{\includegraphics[width=0.4\textwidth]{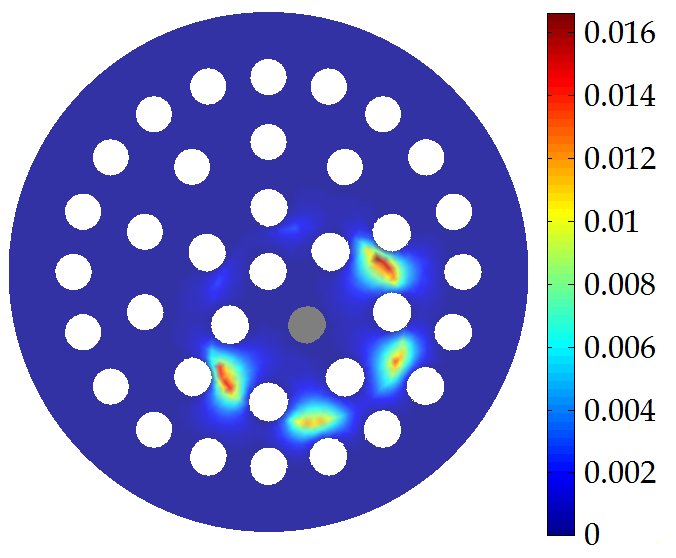}
}
\caption{Global versus localized characteristic functions. Here we use $\delta=0.3$ and obtained maximum error  of $0.016$.}
\label{fig:GlobalCharFunction}
\end{figure}
 To show the effectiveness of the numerical methodology described in section \ref{sec:approximation}, we first consider the problem configuration used in Section \ref{sec:illustration}. See Figure \ref{fig:figure1}. 
We study the expansion which localizes the harmonic characteristic functions. 
We first compare the global and the localized harmonic characteristic functions. 
In Figure \ref{fig:GlobalCharFunction} we plot the global characteristic function (left picture) corresponding to a randomly selected inclusion.
 By construction this harmonic characteristic function is zero at the boundary of all other inclusions which imply 
 a fast decay away from the inclusion. In addition, we plot the difference between the localized 
 characteristic function in
(\ref{eq:deltachi}) and the characteristic function in (\ref{chimulti}); see 
Figure \ref{b2}.  Here we use  $\delta=0.3$ 
and observe that the maximum value of this absolute difference is $0.016$.

\begin{table}[!h]
\begin{center}\footnotesize
\begin{tabular}{|c|c|c|c|} \hline\hline
$\delta$ & $e(u_0-u_{0}^{\delta})$ & 
$e(u_{0,0}-u_{0,0}^{\delta})$&$e(u_{c}-u_{c}^{\delta})$ \\
\hline\hline
 0.001 & 0.830673 & 0.999907& 0.555113\\ 
  0.05 & 0.530459 & 0.768135& 0.549068\\ 
  0.10 & 0.336229 & 0.639191& 0.512751\\ 
  0.20 & 0.081500 & 0.261912& 0.216649\\ 
  0.30 & 0.044613 & 0.088706& 0.048173\\ 
  0.40 & 0.041061 & 0.047743& 0.007886\\ 
  0.50 & 0.033781 & 0.034508& 0.001225\\ 
  0.60 & 0.029269 & 0.029362& 0.000174\\ 
  0.70 & 0.020881 & 0.020888& 0.000021\\ 
  0.80 & 0.012772 & 0.012773& 0.000003\\ 
  0.90 & 0.006172 & 0.006172& 0.000000\\  \hline
  \hline
\end{tabular}\\
\end{center}
\caption[Relative error in the approximation of $u_0$ by using locally computed 
basis functions and truncated boundary conditions. 
$e(w)=\|w\|_{H^1}/\|u_0\|_{H^1}$.]{ Relative error in the approximation of $u_0$ by using locally computed 
basis functions and truncated boundary condition effect. Here 
$u_0=u_{0,0}+u_c$ where $u_c$ is combination of 
harmonic characteristic functions and $u_{0}^{\delta}=u_{0,0}^{\delta}+u_{c}^{\delta}$
is computed by solving $u_{0,0}^{\delta}$ on a $\delta-$strip of the boundary 
$\partial D$ and the basis functions on a $\delta-$strip of the 
boundary of each inclusion.}\label{tab:deltaerror}
\end{table}

\begin{figure}[!htb]
\centering
\includegraphics[width=\textwidth]{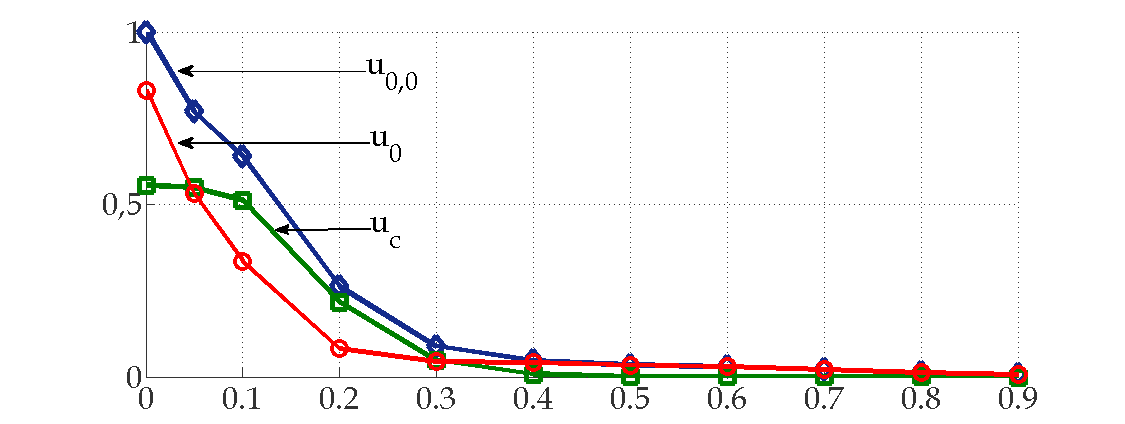}
\caption{Relative error in the approximation of $u_0$ given by Table \ref{tab:deltaerror}}
\label{fig:deltaerror}
\end{figure}

For further comparison, we introduce the relative $H^1$ error from $u_0$ in (\ref{u0multiscale}) 
to its approximation $u_0^\delta$ in (\ref{u0trnmultiscale}). This  is given by 
\begin{eqnarray*}
e(u_0-u_0^{\delta})&=&\frac{\|u_0-u_0^{\delta}\|_{H^1}}{\|u_0\|_{H^1}}.
\end{eqnarray*}
Analogously, the relative $H^1$ error of the approximation of $u_{0,0}$ is given by
\begin{eqnarray*}
e(u_{0,0}-u_{0,0}^{\delta})&=&\frac{\|u_{0,0}-u_{0,0}^{\delta}\|_{H^1}}{\|u_{0,0}\|_{H^1}}.
\end{eqnarray*}
The error $e\left(u_c-u_c^{\delta}\right)$ is defined in a similar way.
According to Theorem \ref{theorem}, the error between  the exact solution of  problem (\ref{Strongform})
 with coefficient (\ref{eq:coeff1inc-multiple}) and $u_0$ in (\ref{u0multiscale}) is of order ${\eta}^{-1}$.

Table \ref{tab:deltaerror} and Figure \ref{fig:deltaerror} show the relative errors of the multiscale method.
In Table \ref{tab:deltaerror}, we observe that as the neighborhood 
size grows, the error is  reduced. For instance, for $\delta = 0.2$ the error between the exact solution $u_0$ and the
 truncated solution $u_0^{\delta}$ is $8\%$. By selecting $\delta=0.3$ we obtain a relative error of the order of 4\%. Numerically we observe that 
an optimal value for $\delta$ is the smallest distance that includes one layer of inclusions away from the selected
 one. Therefore, 
this approximation is more efficient for densely packed inclusions. 

\begin{figure}[!htb]
\centering
\includegraphics[width=0.5\textwidth]{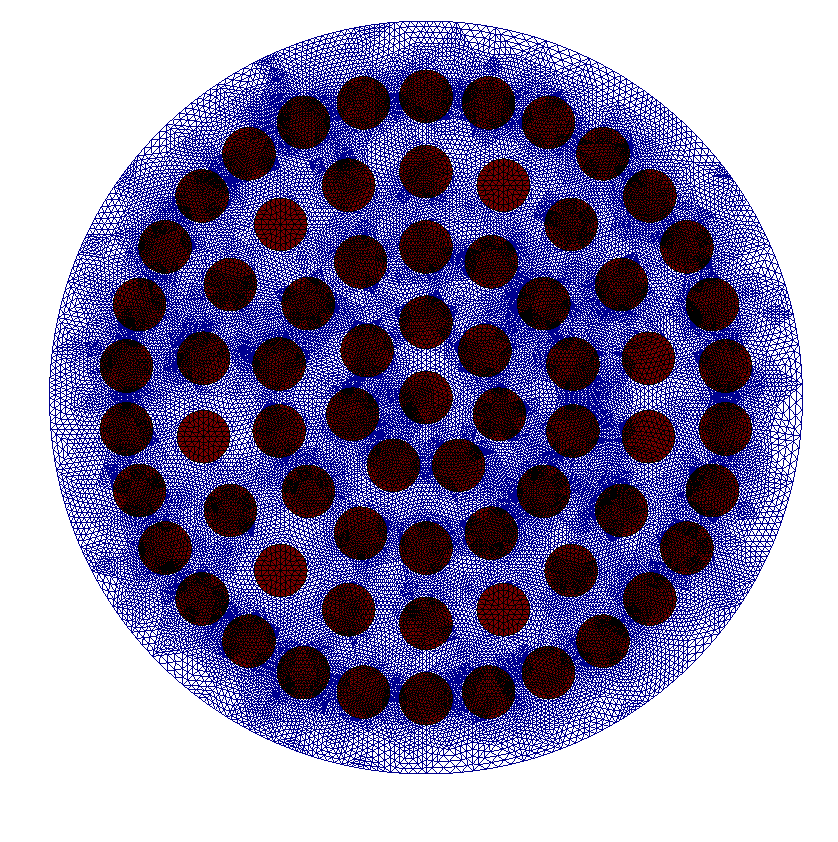}
\caption{Geometry for the problem (\ref{Example2Multiscale}).}
\label{fig:decay} 
\end{figure}

We now consider an additional geometrical configuration of inclusions. 
We consider  $D=(0,1)$,  the circle with center $(0,0)$ and radius $1$, and $60$ 
(identical) circular inclusions of radius $0.07$. 
We consider the problem,
\begin{equation}\label{Example2Multiscale}
\left\{\begin{array}{ll}
-\mbox{div}(\kappa(x)\nabla u_\eta(x))=1, & \mbox{ in }D\\
\hspace{0.1in}u(x)=x_1+x_2^2, & \mbox{ on }\partial D,
\end{array}\right.
\end{equation}
In the Figure \ref{fig:decay} we illustrate the geometry for problem (\ref{Example2Multiscale}). Table \ref{tab:deltaerror2} summarizes similar results to those described above.

\begin{table}[!htb]
\begin{center}\footnotesize
\begin{tabular}{|c|c|c|c|} \hline\hline
$\delta$ & $e(u_0-u_{0}^{\delta})$ & 
$e(u_{00}-u_{0,0}^{\delta})$&$e(u_{c}-u_{c}^{\delta})$ \\
\hline\hline
0.001 & 0.912746 & 0.999972& 0.408063\\ 
  0.05 & 0.369838 & 0.549332& 0.399472\\ 
  0.10 & 0.181871 & 0.351184& 0.258946\\ 
  0.20 & 0.013781 & 0.020172& 0.011061\\ 
  0.30 & 0.013332 & 0.013433& 0.000737\\ 
  0.40 & 0.010394 & 0.010396& 0.000057\\ 
  0.50 & 0.009228 & 0.009228& 0.000004\\ 
  0.60 & 0.006102 & 0.006102& 0.000000\\ 
  0.70 & 0.005561 & 0.005561& 0.000000\\ 
  0.80 & 0.002239 & 0.002239& 0.000000\\ 
  0.90 & 0.001724 & 0.001724& 0.000000\\ \hline
\hline
\end{tabular}\\

\end{center}
\caption{ Relative error in the approximation of $u_0$ when using locally computed 
basis functions and a truncated boundary conditions. 
$e(w)=\|w\|_{H^1}/\|u_0\|_{H^1}$. Here 
$u_0=u_{0,0}+u_c$ where $u_c$ is combination of 
harmonic characteristic functions and $u_{0}^{\delta}=u_{0,0}^{\delta}+u_{c}^{\delta}$
is computed by solving $u_{0,0}^{\delta}$ on a $\delta-$strip of the boundary 
$\partial D$ and the basis functions on a $\delta-$strip of the 
boundary of each inclusion.}\label{tab:deltaerror2}
\end{table}

\begin{figure}[!h]
\centering
\includegraphics[width=\textwidth]{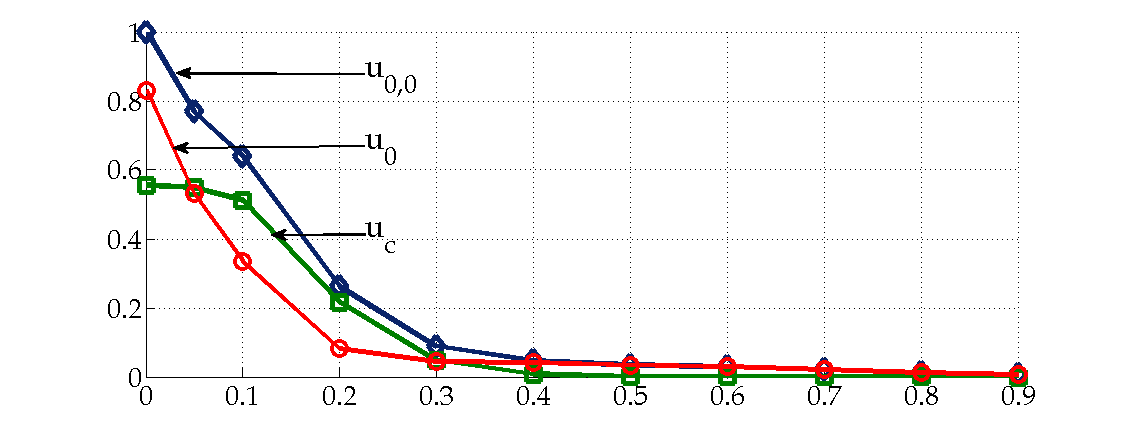}
\caption{Relative error in the approximation of $u_0$ given by Table \ref{tab:deltaerror2}}
\label{fig:deltaerror2}
\end{figure}

The problem setup and simulation of the localized harmonic characteristic functions proceeds
 as described in Section 5. In the present setup the inclusions are clustered more tightly. This induces a 
faster decay of the characteristic 
harmonic function for each individual inclusion, see Figure \ref{fig:deltaerror2}. Thus, as expected in Table \ref{tab:deltaerror2} 
 we observe a reduction in the 
relative error when
 compared to Table \ref{tab:deltaerror}  for a fixed value of $\delta$. For example, $\delta=0.2$
 induces a relative error of $1 \%$ on 
 $u_0^{\delta}$ and of $2 \%$ for $u_{0,0}^{\delta}$ for the geometry show in the Figure \ref{fig:decay} 
while this value of $\delta$  induces errors 
 of $8 \%$ and $26 \%$ for these two variables for the geometry shown in Figure \ref{a1}.

\section{Conclusions}\label{conclusions}

We consider the solution of elliptic problems modeling properties of composite materials. 
Using an expansion in terms
 of the properties ratio presented in \cite{calo2014asymptotic}, we design a 
multiscale method to approximate  solutions. We develop procedures 
that effectively and accurately compute the first few terms in the expansion. In particular, we compute the 
asymptotic limit which is an approximation of order $\eta^{-1}$ to the solution (where $\eta$ represent 
the ratio between lowest and highest material property values). The expansion in \cite{calo2014asymptotic}
is written in terms of the harmonic characteristic functions that are globally supported functions,
 one for each inclusion.
 The main idea we propose is  to approximate the harmonic characteristic functions by solving
 local problems around each inclusion. We use numerical examples to compute the asymptotic
 limit $u_0$ with the localized harmonic characteristic functions.
The analysis of the truncation error depends on decay properties of the harmonic characteristic
 functions and is under current investigation. This method can be used in several important 
engineering applications with heterogeneous coefficients such as complex flow in porous media and complex modern materials.

\bibliographystyle{siam} 
\bibliography{PANBibAMsFEM}

\end{document}